
%
%

%
%

%
%


\font\scriptsize=cmr7

\input amstex

\documentstyle{amsppt}
  \magnification=1000
  \hsize=7.0truein
  \vsize=9.0truein
  \hoffset -0.1truein
  \parindent=2em

\define\Ac{{\Cal A}}                           

\define\Aco{\overset {\scriptsize o}\to \Ac}   

\define\Acw{\widetilde{\Ac}}                   

\define\Ao{\overset {\scriptsize o}\to A}      



\define\aw{\tilde{a}}                          


\define\Bc{{\Cal B}}                           

\define\Bcw{\widetilde{\Bc}}                   

\define\betaw{\tilde{\beta}}                   

\define\Bo{\overset {\scriptsize o}\to B}      

\define\card{\text{\rm card}}                  

\define\Cc{{\Cal C}}                           

\define\Ccw{\widetilde{\Cc}}                   

\define\Co{\overset {\scriptsize o}\to C}      

\define\Cpx{\bold C}                           

\define\crc{^{\scriptsize o}}                  

\define\cw{\tilde{c}}                          

\define\Dc{{\Cal D}}                           


\define\dcup#1{\underset#1\to\cup}             

\define\fdim{\text{\rm fdim}\,}                

\define\freeF{\bold F}                         

\define\gammab{\bar{\gamma}}                   

\define\gammaw{\tilde{\gamma}}                 

\define\Ghat{\widehat G}                       

\define\Indeed{\demo{Indeed}}                  

\define\Intg{\bold Z}                          

\define\isorarrow{\overset\sim\to\rightarrow}  

\define\lspan{\text{\rm span}@,@,@,}           

\define\lw{\tilde{l}}                          

\define\MvN{{\Cal M}}                          

\define\Nats{\bold N}                          

\define\NvN{{\Cal N}}                          



\define\Pc{{\Cal P}}                           

\define\Pf{\demo{Proof}}                       


\define\piot{\tfrac\pi2}                       

\define\pol{\text{\rm pol}}                    

\define\psiw{\widetilde{\psi}}                 

\define\pw{\tilde{p}}                          

\define\QED{$\hfill$\qed\enddemo}              

\define\Qo{\overset {\scriptsize o}\to Q}      


\define\Reals{\bold R}                         

\define\restrict{\lower .3ex                   
     \hbox{\text{$|$}}}

\define\Ro{\overset {\scriptsize o}\to R}      

\define\Rw{\widetilde{R}}                      

\define\smd#1#2{\underset{#2}\to{#1}}          

\define\smdb#1#2{\undersetbrace{#2}\to{#1}}    

\define\smdbp#1#2#3{\overset{#3}\to            
     {\smd{#1}{#2}}}

\define\smdbpb#1#2#3{\oversetbrace{#3}\to      
     {\smdb{#1}{#2}}}

\define\smdp#1#2#3{\overset{#3}\to             
     {\smd{#1}{#2}}}

\define\smdpb#1#2#3{\oversetbrace{#3}\to       
     {\smd{#1}{#2}}}

\define\Sw{\widetilde{S}}                      

\define\SwpbSw{\Sw'\backslash\Sw}              

\define\tf{\text{followed by}}                 

\define\tp{\text{preceded by}}                 

\define\tq{\text{"}}                           

\define\tqc{\hfil\text{"}\hfil}                

\define\uw{\tilde{u}}                          

\define\vw{\tilde{v}}                          

\define\Zw{\widetilde{Z}}                      

\topmatter

  \title Free products of hyperfinite von~Neumann algebras
         and free dimension \endtitle

  \author Ken Dykema \endauthor

  \affil University of California, \\
         Berkeley, California, USA 94720, \\
         (e--mail dykema\@math.berkeley.edu) \endaffil

  \thanks Studies and research
          supported by the Fannie and John Hertz Foundation. \endthanks
  \thanks This work will form part of the author's Ph.D\. thesis
          at the University of California, Berkeley. \endthanks

  \date April 1992 \enddate

\endtopmatter

\document
  \TagsOnRight
  \baselineskip=14pt

\noindent{\bf Introduction.}

  Voiculescu's theory of freeness in noncommutative probability spaces
    (see~\cite{10,11,12,13,14,15}, especially the latter for
    an overview) has made possible the recent surge of results
    about and related to the free group factors
    $L(\freeF_n)$~\cite{13,3,7,8,9,4}.
  One hopes to eventually be able to solve the old isomorphism question,
    first raised by R.V. Kadison in the 1960's,
    of whether
    $L(\freeF_n)\cong L(\freeF_m)$ for $n\neq m$.

  In Voiculescu's theory, (see also~\cite{1}), one takes free products
    of finite von Neumann algebras, denoted $\Ac*\Bc$, and one has
    $L(G)*L(H)\cong L(G*H)$, (where $L(G)$ for $G$ a discrete group is
    the group von Neumann algebra, which is generated by the left regular
    representation of $G$ on $l^2(G)$).
  It is of  intrinsic interest to decide when $\Ac*\Bc$ is a factor, and
    to determine its isomorphism class.
  It may be that such results or the techniques used will give insight into
    the isomorphism problem.
  Moreover, such free products are related to amalgamated free products,
    which have arisen in connection with results about
    irreducible subfactors~\cite{6,9,1$'$}.

  In~\cite{4} and~\cite{9}, the interpolated free group factors
    $L(\freeF_r)$ $(1<r\le\infty)$ were found, that have equality with the
    free group factor on $n$ generators if $r=n\in\Nats\backslash\{0,1\}$
    and that satisfy
    $$ L(\freeF_r)*L(\freeF_{r'})=L(\freeF_{r+r'})
      \qquad(1<r,r'\le\infty) \tag1 $$
    and
    $$ L(\freeF_r)_{\gamma}=L(\freeF(1+\frac{r-1}{\gamma^2}))
      \qquad(1<r\le\infty,\,0<\gamma<\infty). \tag2 $$
  Moreover, it was seen that these two formulas imply that either all
    the $L(\freeF_r)$ are isomorphic to each other or that no two
    of them are isomorphic.

  In this paper, we examine free products of general hyperfinite
    von Neumann algebras
    and express the
    answer in terms of the interpolated free group factors.
  (A finite von Neumann algebras is one which
    has an everywhere defined faithful trace--state;
   it is hyperfinite if it has a dense subalgebra which is the increasing union
    of finite dimensional algebras,
    {\it i.e\.}
    for us ``hyperfinite'' means finite approximately finite dimensional;
   since by~\cite{5} there is only one hyperfinite II$_1$--factor $R$,
    every hyperfinite von Neumann algebra is a direct sum of $R$ tensor
    an abelian algebra and $M_n$ tensor an abelian algebra, as $n$ ranges
    over $\Nats\backslash\{0\}$.)
  Our results are best organized using the concept of
    {\it free dimension}.
  A von Neumann algebra $\Ac$ belonging to a certain class of algebras and
    having specified faithful trace--state
    is assigned a positive real number as its free dimension,
    denoted $\fdim(\Ac)$.
  For example, $\fdim(L(\freeF_r))=r$, $\fdim(L(\Intg))=\fdim(R)=1$,
    where $R$ is the hyperfinite II$_1$~factor.
  Now suppose that $\Ac$ and $\Bc$ are hyperfinite von Neumann
    algebras having specified faithful trace--state.
  We give precise conditions for $\Ac*\Bc$ to be a factor,
    (essentially, when $\Ac$ and $\Bc$ are not too ``lumpy''),
    and show that when
    we have factoriality, $\Ac*\Bc=L(\freeF(\fdim(\Ac)+\fdim(\Bc)))$.
  When we do not have factoriality, we have
    $\Ac*\Bc\cong L(\freeF_r)\oplus C$, where $C$ is a finite
    dimensional algebra, and
    $$ \fdim(\Ac*\Bc)=\fdim(\Ac)+\fdim(\Bc). \tag3 $$

 $C$ can be found by examining the minimal projections of $\Ac$ and $\Bc$,
and $r$ is found using (3).  Note that if the free group factors are
isomorphic to each other, "free dimension" ceases to be well-defined,
but also the need to find $r$ in the above expression disappears.

  The proofs given here
    consist in large part of algebraic manipulations of the sort
    introduced in~\cite{4}, but we also make use of the two pictures
    of $L(\freeF_r)$ contained in~\cite{9} and~\cite{4}.
  Both of these pictures depend on Voiculescu's random matrix model for
    freeness~\cite{14}.

  This paper has five sections.
  In~\S1 we establish notation and prove some theorems that will be
    our basic tools in the sequel;
   in~\S2 we examine free products of finite dimensional abelian algebras
    and related algebras;
   in~\S3 we examine free products of general finite dimensional algebras;
   in~\S4 we develop inductive limit techniques in order to extend the
    previous results to free products of general hyperfinite
    von Neumann algebras;
   in~\S5 we write down what these results imply for the free products of
    amenable groups, using Connes' charactarization in~\cite{2}.

\noindent{\bf \S1. Notation and basic theorems.}

  Throughout this paper, we will be working with finite von Neumann algebras,
    and each one will have a normalized, faithful trace associated to it.
  Moreover, when we write that two von Neumann algebras are isomorphic,
    (or one is contained in the other), we will mean that the
    isomorphism (inclusion) is trace--preserving.
  To save writing, we will not always describe these traces explicitly,
    but will rely on the  conventions found in Remark~3.1
    of~\cite{4}.
  In addition, we will use the following notation to specify the trace on
    a direct sum of algebras:  for von Neumann algebras $A$ and $B$
    with traces $\tau_A$ and $\tau_B$, $\smd A\alpha \oplus \smd B\beta$
    where $\alpha,\beta\ge0$ and $\alpha+\beta=1$, will denote the
    algebra $A\oplus B$ whose associated trace is
    $\tau(a,b)=\alpha\tau_A(a)+\beta\tau_B(b)$.
  In the case that one of the numbers $\alpha$, $\beta$ equals zero,
    say for example $\beta=0$, we will take $\smd A\alpha \oplus \smd B\beta$
    to denote the algebra $A$ with trace $\tau_A$.
  If we also want to give names to the projections corresponding to
    the identity elements of $A$ and $B$, we will write
    $$ \smdp A\alpha p \oplus \smdp B\beta q, $$
    meaning $p=(1,0)$ and $q=(0,1)$.
  Similar notation will apply to direct sums of more that two algebras.

  Traveling (and alternating) products
    will be as defined
    in~3.4 of~\cite{4}, and here also $\Lambda(\ldots)$ will
    denote the set of all traveling products.
  Also, for an algebra $A$, $\Ao$ will denote the kernel of the trace on $A$.

  We will be using the interpolated free group factors~\cite{4,9},
    $L(\freeF_r)$ for $1<r<\infty$, and for convenience
    we will also let $L(\freeF_1)$ denote $L(\Intg)$.

  \proclaim{Theorem 1.1}(Two projections.)
    Let $1>\alpha\ge\max\{\beta,1-\beta\}>0$.
    Then
      $$ (\smdp\Cpx\alpha p\oplus\smdp\Cpx{1-\alpha}{1-p})
        *(\smdp\Cpx\beta q\oplus\smdp\Cpx{1-\beta}{1-q})
        \cong \smdp\Cpx{\alpha+\beta-1}{p\wedge q}\oplus
        \bigl(\undersetbrace{2(1-\alpha)}
        \to{L^\infty([0,\piot ],\nu)\otimes M_2}\bigr)
        \oplus\smdp\Cpx{\alpha-\beta}{p\wedge(1-q)}, \tag4 $$
      where $\nu$ is a probability measure without atoms on $[0,\piot ]$,
      and $L^\infty([0,\piot],\nu)$ has trace given by integration with
      respect to $\nu$.
    In the picture of the right hand side of~(4), we have
      $$ \aligned
        p&=1\oplus\left(\smallmatrix 1&0\\0&0\endsmallmatrix\right)
          \oplus1, \\
        q&=1\oplus\left(\smallmatrix \cos^2\theta&\cos\theta\sin\theta\\
          \cos\theta\sin\theta&\sin^2\theta\endsmallmatrix\right)
          \oplus0,
      \endaligned $$
      where $\theta\in[0,\piot]$.
  \endproclaim
  \Pf
    The proof is the same as that of Proposition~3.2 of~\cite{4}.
    The algebra on the left hand side of~(4) has a dense subalgebra which is
      a quotient of the universal unital C$^*$--algebra generated by
      two projections.
    The trace on this subalgebra is then determined using Voiculescu's
      result~\cite{11} about the distribution of $pqp$.
    (We do the two cases $\beta\ge\frac12$ and $\beta\le\frac12$ separately.)
  \QED

  \proclaim{Theorem 1.2}
    Let $A$, $B$ and $C$ be finite von Neumann algebras (with
      implicitly specified traces) and $n\ge1$.
    Let
      $$ \matrix\format\l&\c&\l\\
        \MvN=\bigl(&\smdbpb{(C\otimes M_n)}\alpha p&\oplus A\bigr)*B \\
        \cup &&\\
        \NvN=\bigl(&\smdp{M_n}\alpha p&\oplus A\bigr)*B \\
      \endmatrix $$
      and let $p_1\le p$ be a minimal projection of $M_n$.
    Then in $p_1\MvN p_1$ we have that $p_1\NvN p_1$
      and $C\otimes p_1$ are free, and together they generate $p_1\MvN p_1$,
      so
      $$ p_1\MvN p_1\cong C*p_1\NvN p_1. $$
    Moreover, the central projection of $p_1$ in $\MvN$ equals the
      central projection of $p_1$ in $\NvN$.
  \endproclaim
  \Pf
    For notational convenience, we identify $C$
      with $(C\otimes1)\oplus0\subset\MvN$.
    To see that $p_1\NvN p_1$ and $p_1 C$ generate $p_1\MvN p_1$,
      note that $\NvN$ and $p_1 C$ generate $\MvN$, so
      $\lspan\Lambda(\NvN,p_1 C)$ is dense in $\MvN$
      and $p_1\Lambda(\NvN,p_1 C)p_1=\Lambda(p_1\NvN p_1,p_1 C)$.

    Now we shall show that a nontrivial alternating product in
      $(p_1\NvN p_1)\crc$ and $p_1\Co$ has trace zero.
    Let $a=p_1-\frac\alpha n$.
    Then $(M_n\oplus A)\crc=\Cpx a+S$
      where $S=\{s\in M_n\oplus A\mid\tau(s)=0,\,p_1sp_1=0\}$.
    So $\lspan\Lambda(\{a\}\cup S,\Bo)$ is a $*$--algebra that is
      dense in $\NvN$.
    Let $x\in(p_1\NvN p_1)\crc$.
    Then by the Kaplansky density theorem
      $x$ is the s.o.--limit of a bounded sequence $\{R_k\}_{k=1}^\infty$
      in
      $\lspan\Lambda(\{a\}\cup S,\Bo)$.
    For
      $Q\in\lspan(\Lambda(\{a\}\cup S,\Bo)\backslash\{1\})$,
      its easy to show that the trace of $p_1Qp_1$ is equal to a fixed constant
      times the coefficient of $a$ in $Q$.
    So since the trace of $R_k$ and the trace of $p_1R_kp_1$ tend to zero,
      we may assume that the coefficients in each $R_k$ of $1$ and $a$
      are zero.
    Since $R_k-p_1R_kp_1$ tends to zero, we may assume that also the
coefficient
      of each element of $S$ in $R_k$ is zero,
      {\it i.e\.} that each
      $R_k\in\lspan(\Lambda(\{a\}\cup S,\Bo)\backslash(\{1,a\}\cup S))$.
    To prove the theorem it thus suffices to show that a nontrivial alternating
      product in
      $\Lambda(\{a\}\cup S,\Bo)\backslash(\{1,a\}\cup S)$ and $p_1\Co$
      has trace zero.
    But regrouping and multiplying some neighboring elements gives
      (a constant times) a nontrivial alternating product in
      $\{a\}\cup S\cup(p_1\Co)$ and $\Bo$, which has trace zero
      by freeness.

    Let $l$ be the central projection of $p_1$ in $\NvN$.
    Then $l\ge p$ and $l$ commutes with both $B$ and $M_n\oplus A$.
    Clearly the central projection of $p_1$ in $\MvN$ contains $l$.
    But as $l$ commutes with both $(C \otimes M_n)\oplus A$ and $B$,
      $l$ is in the center of $\MvN$, so the central projection
      of $p_1$ in $\MvN$ equals $l$.
  \QED

  By using these two theorems and paying attention to central projections,
    one can calculate the free products of direct sums of algebras,
    provided one knows the free products of the summands with various
    algebras.
  This idea is at the heart of much of what follows in this paper,
    but we get started using R\u{a}dulescu's picture of the interpolated
    free group factors~\cite{9}
    and some techniques that he developed.

  \proclaim{Lemma 1.3}
    In a W$^*$--probability space $(\Cc,\psi)$ with $\psi$ a trace,
      let $\{a_1,a_2,c,u\}$ be $*$--free with $a_1$ and $a_2$ semicircular,
      $c$ circular and where $u$ generates a diffuse abelian von Neumann
      algebra.
    Let $l\in\{u\}''$ be a self--adjoint projection having trace $t$,
      and let $v=\pol(lcl)$, a partial isometry from $l$ to $l$, such that
      $lcl=v(lc^*lcl)^{1/2}$.
    Consider $\Dc=\{a_2,u,v^*a_1v,v^*c\}''$.
    Then $\Dc=L(\freeF(2+2t))$.
  \endproclaim

  \Pf
    We will show that, by adjusting $\Cc$ is necessary, we may assume that
      there exist semicircular elements $Z_1$, $Z_2$ and $Z_3$ in $\Cc$
      such that $\{a_2,u,Z_1,Z_2,Z_3\}$ is a $*$--free family,
      $$ \aligned
        & v^*a_1v=lZ_1l, \\
        & v^*c(1-l)=lZ_2(1-l), \\
        & \{v^*cl\}''=\{lZ_3l\}''.
      \endaligned \tag5 $$
    If $t=1/k$ for $k\in\Nats\backslash\{0\}$, we may apply Voiculescu's
      matrix model as Voiculescu did in Theorem~3.3 of~\cite{13}
      to find $Z_1$, $Z_2$ and $Z_3$ as required.
    For general $t$, let $k\in\Nats$ be such that $t\ge1/k$.
    Consider a W$^*$--probability space $(\Ccw,\psiw)$ in which there is
      a $*$--free family $\{\aw_1,\aw_2,\cw,\uw\}$ with $\aw_1$ and $\aw_2$
      semicircular, $\cw$ circular, $\uw$ generating a diffuse abelian
      von Neumann algebra, and let $\lw$, $\pw$ be projections in $\{\uw\}''$
      with $\lw\le\pw$ and $\psiw(\lw)=\frac1k$, $\psiw(\pw)=\frac1{tk}$.
    Let $\vw=\pol(\lw\cw\lw)$.
    We may now assume that there are semicircular $\Zw_1$, $\Zw_2$ and $\Zw_3$
      in $\Ccw$ such that $\{\aw_2,\uw,\Zw_1,\Zw_2,\Zw_3\}$ is $*$--free and
      $$ \aligned
        & \vw^*\aw_1\vw=\lw\Zw_1\lw, \\
        & \vw^*\cw(1-\lw)=\lw\Zw_2(1-\lw), \\
        & \{\vw^*\cw\lw\}''=\{\lw\Zw_3\lw\}''.
      \endaligned $$
    Let $a_i=\pw\aw_i\pw$ $(i=1,2)$, $c=\pw\cw\pw$, $u=\pw\uw$.
    By Proposition~2.3 of~\cite{13}, $\{a_1,a_2,c,u\}$ is $*$--free in
      $(\pw\Ccw\pw,tk\psiw\restrict_{\pw\Ccw\pw})$.
    Letting $Z_i=\pw\Zw_i\pw$ $(i=1,2,3)$ we have also by the same proposition
      that in $\pw\Ccw\pw$ each $Z_i$ is semicircular and
      $\{a_2,u,Z_1,Z_2,Z_3\}$ is $*$--free.
    Writing $\lw=l$ (which in $\pw\Ccw\pw$ has trace $t$),
      clearly~(5) holds, as required.

    So $\Dc\cong\{a_2,u,lZ_1l,lZ_2(1-l),lZ_3l\}''$, which by R\u{a}dulescu's
      definition of the interpolated free group factors~\cite{9}
      is $L(\freeF(2+2t))$.
  \QED

  \proclaim{Lemma 1.4}
    For $0\le \delta\le1$ we have
      $$ L(\Intg)*(\smd{L(\Intg)}\delta\oplus\smd\Cpx{1-\delta})
        =L(\freeF(1+\delta^2+2\delta(1-\delta))). \tag6 $$
  \endproclaim
  \Pf
    Let us first prove the case $\delta\ge1/2$.
    Let $\Ac$ denote the algebra of the left hand side of~(6),
      and let $\tau$ be its trace.
    We can model $\Ac$ as follows.
    Let $(\MvN,\tau)=(\NvN\otimes M_2,\phi\otimes\tau_2)$,
      where $(\NvN,\phi)$ is a W$^*$--probability space containing enough
      free elements (to be specified below), $\phi$ is a trace
      and $\tau_2$ is the normalized trace on $M_2$.
    In $\MvN$
      let $X=\left(\smallmatrix a_1&c\\c^*&a_2\endsmallmatrix\right)$,
      $Y=\left(\smallmatrix 0&0\\0&u\endsmallmatrix\right)$,
      where in $(\NvN,\phi)$, $a_1$ and $a_2$ are semicircular,
      $c$ is circular, u is a Haar unitary
      and $\{a_1,a_2,c,u\}$ is $*$--free.
    Let $l\in\{u\}''$ in $\NvN$ be a self--adjoint projection of trace
      $\frac1\delta-1$.
    Let $p=\left(\smallmatrix l&0\\0&1\endsmallmatrix\right)$,
      so $\tau(p)=(2\delta)^{-1}$.
    In $(p\MvN p,2\delta\tau\restrict_{p\MvN p})$ consider the
      elements $pXp$ and $pYp$.
    By~2.3 of~\cite{13}, they are $*$--free and $pXp$ is
      semicircular.
    Clearly $pYp$ generates an abelian von Neumann algebra
      isomorphic to $\smd{L(\Intg)}\delta\oplus\smd\Cpx{1-\delta}$,
      so $\Ac$ is isomorphic to $\Bc=\{pXp,pYp\}''$ in $p\MvN p$.
    Now let $f=\left(\smallmatrix 0&0\\0&1\endsmallmatrix\right)$,
      $e=\left(\smallmatrix 0&0\\0&l\endsmallmatrix\right)$,
      so $f$ has trace $\delta$ and $e$ has trace $1-\delta$ in $p\MvN p$.
    Let $w=\pol((p-f)Xe)=
      \left(\smallmatrix 0&\pol(lcl)\\0&0\endsmallmatrix\right)=
      \left(\smallmatrix 0&v\\0&0\endsmallmatrix\right)$,
      (where $v=\pol(lcl)$), a partial isometry from $e$ to $p-f$.
    Then $f\Bc f=\{fXf,fYf,w^*Xw,w^*Xf\}''$, which is clearly isomorphic
      to the von Neumann algebra $\{a_2,u,v^*a_1v,v^*c\}''$ in $\NvN$.
    Lemma~1.3 gives us that $\Ac_\delta\cong f\Bc f=\L(\freeF(2/\delta))$,
      and formula~(2) gives $\Ac=L(\freeF(1+\delta^2+2\delta(1-\delta)))$
      as required.

    Let us prove the general case by induction on $k$,
      where $\frac1{k+1}\le\delta\le\frac1k$.
    We just proved the case $k=1$.
    For $k\ge2$, let
      $$ \matrix\format\c&\l&\c&\c&\c&\c&\c&\c&\c&\l\\
        \MvN&=(&L(\Intg)&\oplus&L(\Intg)&)*(&L(\Intg)&\oplus&\Cpx&) \\
        \cup \\
        \NvN_2&=(&L(\Intg)&\oplus&L(\Intg)&)*(&\Cpx&\oplus&\Cpx&) \\
        \cup \\
        \NvN_1&=(&\Cpx&\oplus&L(\Intg)&)*(&\Cpx&\oplus&\Cpx&) \\
        \cup \\
        \NvN_0&=(&\smdp\Cpx\delta{1-p}&\oplus&\smdp\Cpx{1-\delta}p&)
          *(&\smdp\Cpx\delta{1-p}&\oplus&\smdp\Cpx{1-\delta}q&).
      \endmatrix $$
    Now by Theorem~1.1,
      $\NvN_0=\smd\Cpx{1-2\delta}\oplus\smd{(L(\Intg)\otimes M_2)}{2\delta}$
      and $p\NvN_0p=\smd\Cpx{\frac{1-2\delta}{1-\delta}}
      \oplus\smd{L(\Intg)}{\frac\delta{1-\delta}}$.
    Since $\frac1k\le\frac\delta{1-\delta}\le\frac1{k-1}$,
      we may use Theorem~1.2 and inductive hypothesis to conclude that
      $p\NvN_1p=L(\freeF(1+(1+\delta)^2(\delta^2+2\delta(1-2\delta))))$.
    Since the central projection of $p$ in $\NvN_0$ is $1$, we conclude
      by~1.2 and~(2) that
      $\NvN_1=L(\freeF(1+2\delta-3\delta^2))$.
    Applying Theorem~1.2 twice again, we get
      $\NvN_2=L(\freeF(1+2\delta-2\delta^2))$ and
      $\MvN=L(\freeF(1+2\delta-\delta^2))$.
  \QED

  \proclaim{Remark 1.5} \rm
    We now have another description of the interpolated free group
      factors $L(\freeF_r)$ ($r\ge1$) by writing
      $$ L(\freeF_r)=L(\freeF_{[r]})
        *(\smd{L(\Intg)}\delta\oplus\smd\Cpx{1-\delta}), $$
      where $[r]$ is the integer part of $r$ and
      $\delta^2+2\delta(1-\delta)=r-[r]$, $0\le\delta<1$.
    Thus $L(\freeF_r)=\L(\Intg)*D_{r-1}$ where $D_{r-1}$ equals
      $L(\freeF_{r-1})$ if $r\ge2$ and $\smd{L(\Intg)}\delta\oplus\Cpx$
      if $r<2$.
  \endproclaim

  \proclaim{Lemma 1.6}
    For $s\ge1$ and $0\le\delta\le1$ we have
      $$ L(\Intg)*(\smd{L(\freeF_s)}\delta\oplus\smd\Cpx{1-\delta})
        =L(\freeF(1+s\delta^2+2\delta(1-\delta))). \tag7 $$
  \endproclaim
  \Pf
    First show that $L(\Intg)*(\smd\Cpx\alpha\oplus\Cpx)=
      L(\freeF(1+2\alpha(1-\alpha)))$ by taking $\alpha\ge1/2$ and writing
      $$ \matrix\format\c&\l&\c&\c&\c&\c&\c&\c&\c&\l\\
        \MvN&=(&\smd{L(\Intg)}\alpha&\oplus&L(\Intg)&)
          *(&\smd\Cpx\alpha&\oplus&\Cpx&) \\
        \cup \\
        \NvN_1&=(&\Cpx&\oplus&L(\Intg)&)*(&\Cpx&\oplus&\Cpx&) \\
        \cup \\
        \NvN_0&=(&\smdp\Cpx{}{p}&\oplus&\Cpx&)
          *(&\smdp\Cpx{}{q}&\oplus&\Cpx&)
          \,=\,\smdp\Cpx{2\alpha-1}{p\wedge q}\oplus (L(\Intg)\otimes M_2).
      \endmatrix $$
    Now apply~1.2 twice, using~1.4 to find $\NvN_1$.
    One then finds~(7) by applying~1.2 again.
  \QED

  \proclaim{Proposition 1.7}
    For $r,s\ge1$ and $0\le\gamma,\delta\le1$ we have
      $$ \align
        \text{(i) }&\;
          L(\freeF_r)*(\smd{L(\freeF_s)}\delta\oplus\smd\Cpx{1-\delta})
          =L(\freeF(r+s\delta^2+2\delta(1-\delta))) \\ \vspace{1\jot}
        \text{and (ii) }&\;
          (\smd{L(\freeF_r)}\gamma\oplus\smd\Cpx{1-\gamma})
          *(\smd{L(\freeF_s)}\delta\oplus\smd\Cpx{1-\delta})
          =\cases L(\freeF(r\gamma^2+2\gamma(1-\gamma)
            +s\delta^2+2\delta(1-\delta)))&\text{if }\;\gamma+\delta\ge1
            \\ \vspace{1\jot}
          \smdb{L(\freeF((\gamma+\delta)^{-2}(r\gamma^2+s\delta^2
            +4\gamma\delta)))}{\gamma+\delta}\oplus\smd\Cpx{1-\gamma-\delta}
            &\text{if }\;\gamma+\delta\le1. \endcases
      \endalign $$
  \endproclaim
  \Pf
    1.6 gives~(i) for the case $r=1$.
    Now we prove~(ii) for the case $r=s=1$.
    Take
      $$ \matrix\format\c&\l&\c&\c&\c&\l&\c&\c&\c&\l\\
        \MvN&=(&L(\Intg)&\oplus&\Cpx&)*(&L(\Intg)&\oplus&\Cpx&) \\
        \cup \\
        \NvN_1&=(&\Cpx&\oplus&\Cpx&)*(&L(\Intg)&\oplus&\Cpx&) \\
        \cup \\
        \NvN_0&=(&\smd\Cpx\gamma&\oplus&\smd\Cpx{1-\gamma}&)
          *(&\smd\Cpx\delta&\oplus&\smd\Cpx{1-\delta}&),
      \endmatrix $$
      find $\NvN_0$ using Theorem~1.1 and apply Theorem~1.2
      twice (as we did in~1.4) to find $\NvN_1$, then $\MvN$.
    The cases $\gamma\ge\max(\delta,1-\delta)$
      and $\gamma\le\max(\delta,1-\delta)$ are done separately, each
      making use of~1.6.

    Now we can prove~(i) for all $r,s\ge1$.
    Let $L(\freeF_r)=L(\freeF_{[r]})*D_{r-1}$, where $D_{r-1}$
      is as in Remark~1.5.
    Then $L(\freeF_r)*(\smd{L(\freeF_s)}\delta\oplus\smd\Cpx{1-\delta})
      =D_{r-1}*L(\freeF(1+s\delta^2+2\delta(1-\delta))$ by~1.6.
    If $r\ge2$ we are done, and if $s\delta^2+2\delta(1-\delta)\ge1$
      we may apply the description of~1.5 and Lemma~1.4 to yield the
      required result.
    Otherwise $D_{r-1}=\smd{L(\Intg)}{\alpha}\oplus\smd\Cpx{1-\alpha}$
      and $L(\freeF(1+s\delta^2+2\delta(1-\delta)))
      =L(\Intg)*(\smd{L(\Intg)}\beta\oplus\smd\Cpx{1-\beta})$
      for some $\alpha$ and $\beta$.
    Now we apply~(ii) for the case we proved in the preceding paragraph.

    Having~(i) for all $r,s\ge1$, we may now prove~(ii) for all $r,s\ge1$
      using the same argument we did earlier for the case $r=s=1$.
  \QED

  \proclaim{Remark 1.8}\rm
    Knowing~1.5 and~1.7, it is now easy to show that for any $A$,
      if $L(\Intg)*A=L(\freeF_{1+s})$ for some $s\ge0$,
      then $L(\freeF_r)*A=L(\freeF_{r+s})$ for all $r\ge1$.
  \endproclaim

\noindent{\bf \S2. Finite dimensional abelian algebras.}

  \proclaim{Definition 2.1} \rm
    The {\it free dimension} of the algebra
      $$ A=\smd\Cpx{\alpha_1}\oplus\cdots\oplus\smd\Cpx{\alpha_n}
        \;\;\;(n\ge1) $$
      is
      $$ \fdim(A)=\sum\Sb1\le i,j\le n\\i\neq j\endSb \alpha_i\alpha_j. $$
    This includes the cases $\fdim(\Cpx)=0$ and
      $\fdim(\smd\Cpx\alpha\oplus\smd\Cpx{1-\alpha})=2\alpha(1-\alpha)$.
    We also define the free dimension of
      $$ B=\smd{L(\freeF_s)}{\alpha_0}\oplus
        \smd\Cpx{\alpha_1}\oplus\cdots\oplus\smd\Cpx{\alpha_n}
        \;\;\;(s\ge1,\,n\ge0) \tag8 $$
      to be
      $$ \fdim(B)=s\alpha_0^2
        +\sum\Sb0\le i,j\le n\\i\neq j\endSb \alpha_i\alpha_j. $$
    This includes the special case $\fdim(L(\freeF_s))=s$.
  \endproclaim

  Of course if the free group factors are isomorphic to each other then
    the free dimension of the algebra of~(8) can take on a continuum of
    values, hence is non--unique.
  We will sometimes want to let $s$ in~(8) be determined by the requirement
    that $B$ have a certain free dimension, but in the case of non--uniqueness
    of free dimension this is no problem since then all values of $s$ would
    give the same algebra.

  In this section we will show the addition formula~(3) for free products of
    algebras as in~2.1 and also decide in which cases we have a factor.
  In particular, for finite dimensional abelian algebras $A$ and $B$
    we have that $A*B$ is a factor if and only if $\dim A,\,\dim B\ge2$,
    $\dim A+\dim B\ge5$ and the largest mass of an atom in $A$ plus the
    largest mass of an atom in $B$ is $\le1$.
  Our proofs will consist of induction arguments using the results of the
    last section.
  The details are fairly tedious and elementary, so we will only
    give outlines of the proofs.

  \proclaim{Lemma 2.2}
    Let
      $$ A=\smd{L(\freeF_r)}{\alpha_0}\oplus\smdp\Cpx{\alpha_1}{p_1}
        \oplus\cdots\oplus\smdp\Cpx{\alpha_n}{p_n}, $$
      where $n\ge0$, $\alpha_0>0$, $\alpha_1\ge\alpha_2\ge\cdots\alpha_n>0$
      and $r\ge1$
    Let $\frac12\le\beta\le1$.
    Then
      $$ A*(\smdp\Cpx\beta q\oplus\smdp\Cpx{1-\beta}{1-q})\cong
        \cases \smd{L(\freeF_s)}{\gamma_0}
          \oplus\smdp\Cpx{\gamma_1}{p_1\wedge q}
          \oplus\cdots\oplus\smdp\Cpx{\gamma_n}{p_n\wedge q}
          &\text{if }\;\beta\ge\alpha_1 \\
        \smd{L(\freeF_s)}{\gamma_0}\oplus
          \smdp\Cpx{\alpha_1+\beta-1}{p_1\wedge q}
          \oplus\smdp\Cpx{\alpha_1-\beta}{p_1\wedge(1-q)}
          &\text{if }\;\beta\le\alpha_1, \\
        \endcases \tag9 $$
      where $\gamma_i=\max(\beta+\alpha_i-1,\,0)$ for $1\le i\le n$
      and where $s$ is such that the free dimension of the right hand side
      of~(9) equals $\fdim(A)+2\beta(1-\beta)$.
  \endproclaim
  \Pf
    For $n=0$, by Remark~1.8 we need only do the case $r=1$.
    But this was done in the proof of~1.6.
    The case n=1 is similarly proved using 1.2 and 1.7.
    One then proceeds by induction on $n$, proving the inductive step
      by writing
      $$ \matrix\format\c&\l&\c&\l\\
        \MvN&=(\smd{L(\freeF_r)}{\alpha_0}\oplus&\smd\Cpx{\alpha_1}\oplus
          \smd\Cpx{\alpha_2}&\oplus\smd\Cpx{\alpha_3}\oplus\cdots\oplus
          \smd\Cpx{\alpha_n})*(\smd\Cpx\beta\oplus\smd\Cpx{1-\beta}) \\
        \cup \\
        \NvN&=(\smd{L(\freeF_r)}{\alpha_0}\oplus&\smd\Cpx{\alpha_1+
          \alpha_2}&\oplus\smd\Cpx{\alpha_3}\oplus\cdots\oplus
          \smd\Cpx{\alpha_n})*(\smd\Cpx\beta\oplus\smd\Cpx{1-\beta}) \\
      \endmatrix $$
      and using~1.2 and~1.7.
    Note that the case $n=2$ and $\alpha_1+\alpha_2\ge\beta$
      requires an additional induction argument on $k$ such that
      $(k-1)\alpha_0\le\alpha_2\le k\alpha_0$.
  \QED

  \proclaim{Theorem 2.3}
    Let
      $$ A=\smdp\Cpx{\alpha_1}{p_1}\oplus\cdots\oplus\smdp\Cpx{\alpha_n}{p_n}
        ,\;\,
        B=\smdp\Cpx{\beta_1}{q_1}\oplus\cdots\oplus\smdp\Cpx{\beta_m}{q_m} $$
      be finite dimensional abelian algebras, where each $\alpha_i$
      and $\beta_j>0$, $n,m\ge2$ and $n+m\ge5$.
    Then
      $$ A*B\cong L(\freeF_s)\oplus\bigoplus\Sb1\le i\le n\\1\le j\le m\endSb
        \smdp\Cpx{\gamma_{ij}}{p_i \wedge q_j}, $$
      where $\gamma_{ij}=\max(\alpha_i+\beta_j-1,0)$
      and where $s$ is chosen so that $\fdim(A*B)=\fdim(A)+\fdim(B)$.
  \endproclaim
  \Pf
    Proceed by induction on $n+m$.
    Assume $n\ge3$ and prove the initial case $n+m=5$ as well as the inductive
      step by writing
      $$ \matrix\format\c&\l&\c&\l\\
        \MvN&=(&\smd\Cpx{\alpha_1}\oplus
          \smd\Cpx{\alpha_2}&\oplus\smd\Cpx{\alpha_3}\oplus\cdots\oplus
          \smd\Cpx{\alpha_n})
          *(\smd\Cpx{\beta_1}\oplus\cdots\oplus\smd\Cpx{\beta_m}) \\
        \cup \\
        \NvN&=(&\smd\Cpx{\alpha_1+
          \alpha_2}&\oplus\smd\Cpx{\alpha_3}\oplus\cdots\oplus
          \smd\Cpx{\alpha_n})
          *(\smd\Cpx{\beta_1}\oplus\cdots\oplus\smd\Cpx{\beta_m}), \\
      \endmatrix $$
      and applying~1.2 using~2.2.
  \QED

  It will be useful to have also the following proposition, which
    is straightforwardly proved from the above results.
  \proclaim{Proposition 2.4}
    Let
      $$ \aligned
        & A=\smd{L(\freeF_r)}{\alpha_0}\oplus
          \smdp\Cpx{\alpha_1}{p_1}\oplus\cdots\oplus\smdp\Cpx{\alpha_n}{p_n}
          \;\;(n\ge0,\,r\ge1,\,\alpha_0\ge0) \\
        & B=\smd{L(\freeF_s)}{\beta}\oplus
          \smdp\Cpx{\beta_1}{q_1}\oplus\cdots\oplus\smdp\Cpx{\beta_m}{q_m}
          \;\;(m\ge0,\,s\ge1,\,\beta_0\ge0),
      \endaligned $$
      where $\alpha_0+\beta_0>0$.
    Then
      $$ A*B\cong L(\freeF_t)\oplus\bigoplus\Sb1\le i\le n\\1\le j\le m\endSb
        \smdp\Cpx{\gamma_{ij}}{p_i \wedge q_j}, $$
      where $\gamma_{ij}=\max(\alpha_i+\beta_j-1,0)$
      and where $t$ is chosen so that $\fdim(A*B)=\fdim(A)+\fdim(B)$.
  \endproclaim

\noindent{\bf \S3. Finite dimensional algebras.}

  In this section we determine the free product of any two finite dimensional
    algebras.  $\Intg_n$ will denote the cyclic group of order $n$, $M_n$ the
    $n\times n$ complex matrices and for a (finite) factor $\MvN$ and
$0<t\le1$,
    $\MvN_t$ denotes an algebra isomorphic to $p\MvN p$, where $p\in\MvN$ is
    a self--adjoint projection of trace $t$.

  \proclaim{Proposition 3.1}
    Let $\Ac$ be a finite von Neumann algebra such that
      $\NvN=L(\Intg_n)*\Ac$ is a factor.
    Then $\MvN=M_n*\Ac$ is a factor and
      $\MvN_{\frac1n}\cong(\NvN_{\frac1n})*L(\freeF_{n-1})$.
  \endproclaim
  \Pf
    Let $p_1,\ldots,p_n$ be an orthogonal family of minimal projections
      of $M_n$, and let $u=\sum_{1\le j\le n}e^{2\pi ij/n}p_j$.
    Then $\lspan\Lambda(\{u,u^2,\ldots,u^{n-1}\},\Aco)$ is a dense
      $*$--subalgebra of $\NvN$.
    Since $\NvN$ is a factor, there are partial isometries $x_k\in\NvN$
      from $p_k$ to $p_1$, {\it i.e\.} such that
      $x_k^*x_k=p_k$, $x_kx_k^*=p_1$ $(2\le k\le n)$.
    Let $v_k\in M_n$ $(2\le k\le n)$ be partial isometries from $p_k$
      to $p_1$.
    Then $p_1\MvN p_1$ is generated by $p_1\NvN p_1$ together with
      $\{v_kx_k^*\mid2\le k\le n\}$.

    \proclaim{Claim 3.1a}
      Each $v_kx_k^*$ is a Haar unitary in $p_1\MvN p_1$.
    \endproclaim
    \Indeed
      It is clearly a unitary, we need only show that the trace of a positive
        power is zero.
      We show that a nontrivial alternating product in $v_k$ and $x_k^*$
        has trace zero.
      Now since the trace of $x_k^*$ is zero and $x_k^* u=e^{2\pi i/n}x_k^*$,
        we may approximate $x_k^*$ in strong--operator topology
        by a bounded sequence in
        $\lspan(\Lambda(\{u,u^2,\ldots,u^{n-1}\},\Aco)
        \backslash\{1,u,u^2,\ldots,u^{n-1}\})$.
      Hence it suffices to show that a nontrivial alternating product in
        $\{v_k\}$ and
        $\Lambda(\{u,u^2,\ldots,u^{n-1}\},\Aco)
        \backslash\{1,u,u^2,\ldots,u^{n-1}\}$ has trace zero.
      But regrouping and multiplying some neighbors gives (a constant times)
        a nontrivial alternating product in $\{v_k,u,u^2,\ldots,u^{n-1}\}$
        and $\Aco$, which by freeness has trace zero.
      This proves the claim.
    \enddemo

    \proclaim{Claim 3.1b}
      $\{p_1\NvN p_1,\,(v_kx_k^*)_{2\le k\le n}\}$ is $*$--free.
    \endproclaim
    \Indeed
      We need only show that a nontrivial traveling product in
        $(p_1\NvN p_1)\crc$,
        $\bigl\{(v_2x_2^*)^m\mid m\in\Intg\backslash\{0\}\bigr\}$, $\ldots$,
        $\bigl\{(v_nx_n^*)^m\mid m\in\Intg\backslash\{0\}\bigr\}$
        has trace zero.
      Regrouping gives a nontrivial alternating product in
        $$ (p_1\NvN p_1)\crc\cup\bigcup_{2\le k\le n}
          \bigl((p_1\NvN p_1)x_k\cup x_k^*(p_1\NvN p_1)\bigr)\cup
          \bigcup\Sb 2\le i,j\le n \\ i\neq j\endSb x_i^*(p_1\NvN p_1)x_j
        \tag10 $$
        and
        $$ \{v_k,\,v_k^*\mid2\le k\le n\}
          \cup\{v_i^*v_j\mid2\le i,j\le n,\,i\neq j\}. \tag11 $$
      Just as above, we may use the Kaplansky density theorem to show that
        each element of the set~(10) is the limit in strong--operator
        topology of a bounded sequence in
        $\lspan(\Lambda(\{u,u^2,\ldots,u^{n-1}\},\Aco)
        \backslash\{1,u,u^2,\ldots,u^{n-1}\})$.
      So it suffices to show that a nontrivial alternating product in
        $\Lambda(\{u,u^2,\ldots,u^{n-1}\},\Aco)
        \backslash\{1,u,u^2,\ldots,u^{n-1}\}$
        and the set~(11) has trace zero.
      Regrouping and multiplying some neighbors gives (a constant times)
        a nontrivial traveling product in
        $\{v_j,v_j^*\mid2\le j\le n\}\cup\{u,u^2,\ldots,u^{n-1}\}$
        and $\Aco$, which by freeness has trace zero.
      This proves the claim and hence the proposition.
    \enddemo
  \QED

  The above proposition, together with Theorem~1.2 and the results of~\S2
    allows us to determine the free product of finite dimensional algebras
    $A$ and $B$, provided that the largest trace of a minimal projection
    of $A$ plus the largest trace of a minimal projection of $B$
    is $\le1$.
  However, we still cannot say what, for example,
    $M_2*(\smd\Cpx{2/3}\oplus\smd\Cpx{1/3})$ is.
  Hence the following:

  \proclaim{Proposition 3.2}
    Let $n\ge2$, $\frac12\le\alpha<1$.
    Then
      $$ (\smdp\Cpx\alpha p\oplus\smdp\Cpx{1-\alpha}{1-p})*M_n\cong
        \cases L(\freeF(1-n^{-2}+2\alpha(1-\alpha)))
          &\text{if }\;\alpha\le1-n^{-2} \\
        \smdp{M_n}{n^2(\alpha+n^{-2}-1}{p'}
        \oplus\smd{L(\freeF(1+n^{-2}-2n^{-4}))}
          {n^2(1-\alpha)}
          &\text{if }\;\alpha\ge1-n^{-2} \endcases \tag12 $$
      and $p'\le p$.
  \endproclaim
  \Pf
    Let $\MvN$ be the left hand side of~(12).
    If $\alpha\le1-n^{-1}$ then the result follows immediately from
      Theorem~2.3 and Proposition~3.1, so suppose $\alpha>1-n^{-1}$.
    Let $(q_i)_{1\le i\le n}$ be an orthogonal family of minimal
      projections of $M_n$ and let
      $\NvN_1=(\{p\}\cup\{q_i\mid1\le i\le n\})''$.
    Also let $\{q_i\mid1\le i\le n\}\cup\{v_{ij}\mid1\le i,j\le n,\,i\neq j\}$
      be a system of matrix units in $M_n$, where $v_{ij}$ is a partial
      isometry from $q_j$ to $q_i$.
    By Theorem~2.3,
      $$ \NvN_1\cong\cases\smdp\Cpx{\alpha+n^{-1}-1}{p\wedge q_1}\oplus
        \cdots\oplus\smdp\Cpx{\alpha+n^{-1}-1}{p\wedge q_n}\oplus
        \smdbp{L(\freeF(1+\frac{n-2}{n^2}))}{n(1-\alpha)}h&\text{if }\;n\ge3 \\
        \smd\Cpx{\alpha-\frac12}\oplus\smd\Cpx{\alpha-\frac12}
        \oplus\smd{(L(\Intg)\otimes M_2)}{2(1-\alpha)}&\text{if }\;n=2.
        \endcases $$
    Then $q_i=p\wedge q_i+h_i$ for $h_i\le h$ and $h=h_1+\cdots+h_n$,
      $\tau(h_i)=1-\alpha$.

  $\MvN$ is generated by $\NvN_1$ together with $\{v_{1j}\mid 2\le j\le n \}$,
so $q_1\MvN q_1$ is generated by $\{v_{1i}h\NvN_1hv_{j1}\mid1\le i,j\le n\}$
together with $\{v_{1j}(p\wedge q_j)v_{j1}\mid 2\le j\le n\}$.
Note that $v_{1j}(p\wedge q_j)v_{j1}=q_1-v_{1j}h_jv_{j1}$.
Let $y_{1j}\in h\NvN_1h$ be a partial isometry from $h_1$ to $v_{1j}h_jv_{j1}$.
Then $v_{1j}y_{j1}$ is a partial isometry from $h_1$ to $v_{1j}h_jv_{j1}$
and $q_1\MvN q_1$ is generated by $h_1\NvN_1h_1$ together with
$\{v_{1j}y_{j1}\mid2\le j\le n\}$.

    \proclaim{Claim 3.2a}
      $\bigl\{q_1\NvN_1q_1,(\{v_{1j}h_jv_{j1}\})_{2\le j\le n}\bigr\}$
        is free in $q_1\MvN q_1$.
    \endproclaim
    \Indeed
      Let $c_j=v_{1j}(h_j-n(1-\alpha)q_j)v_{j1}$.
      We need only show that a nontrivial traveling product in
        $(q_1\NvN_1q_1)\crc$, $(\{c_j\})_{2\le j\le n}$ has trace zero.
      Let $a=p-\alpha$, $b_j=q_j-n^{-1}$ $(1\le j\le n)$.
      Then
        $\lspan\Lambda(\{a\},\{b_j\mid1\le j\le n\})$
        is a dense $*$--subalgebra of $\NvN_1$.
      If $x\in\NvN_1$, $\tau(x)=0$ and $q_jx=x$ for some $1\le j\le n$,
        then, as one can easily show, $x$ is the s.o.--limit of a
        bounded sequence in
        $\lspan(\Lambda(\{a\},\{b_j\mid1\le j\le n\})
        \backslash\{1,b_1,\ldots,b_n\})$.
      So in particular, each $x\in(q_1\NvN_1q_1)\crc$ and each
        $h_j-n(1-\alpha)q_j$ $(2\le j\le n)$ is the limit of such a sequence.
      Thus we must only show that a nontrivial alternating product
        in
        $\Lambda(\{a\},\{b_j\mid1\le j\le n\})
        \backslash\{1,b_1,\ldots,b_n\}$ and
        $\{v_{1j},v_{j1}\mid2\le j\le n\}$ has trace zero.
      Regrouping gives a nontrivial alternating product in $\{a\}$ and
        $\{b_1,\ldots,b_n\}\cup\{v_{1j},v_{j1}\mid2\le j\le n\}$,
        which by freeness has trace zero.
      This proves Claim~3.2a.

    Now $q_1\MvN q_1$ is generated by $\NvN_2$ together with
      $\{v_{1j}y_{j1}\mid2\le j\le n\}$, where we let
      $$ \NvN_2=\bigl(q_1\NvN_1q_1
        \cup\{v_{1j}h_jv_{j1}\mid2\le j\le n\}\bigr)'' $$
      (with $q_1$ for identity element).
    But
      $$ q_1\NvN_1q_1=\smdp\Cpx{1+n\alpha-n}{p_1\wedge q}\oplus
        \smdp{L(\freeF_{n-1})}{n(1-\alpha)}{h_1}, $$
      so by the above claim we have
      $$ \NvN_2\cong(\smdp{L(\freeF_{n-1})}{n(1-\alpha)}{h_1}\oplus
        \smdp\Cpx{}{p_1\wedge q})
        *(\smdp\Cpx{n(1-\alpha)}{v_{12}h_2v_{21}}\oplus
        \smdp\Cpx{}{v_{12}(q_2\wedge p)v_{21}})
        *\cdots
        *(\smdp\Cpx{n(1-\alpha)}{v_{1n}h_nv_{n1}}\oplus
        \smdp\Cpx{}{v_{1n}(q_n\wedge p)v_{n1}}). $$
    Applying~2.2 gives
      \proclaim{Case 1} $\alpha\le1-n^{-2}$, \endproclaim
      \demo{} then
        $$ \NvN_2=L(\freeF(n^2(1-n+2n\alpha-(1+n)\alpha^2))), $$
      \enddemo
      \proclaim{Case 2} $\alpha\ge1-n^{-2}$, \endproclaim
      \demo{} then
        $$ \NvN_2=\smdp{L(\freeF(2-n^{-1}-n^{-2}))}{n^2(1-\alpha)}k
          \oplus\smdp\Cpx{1-n^2(1-\alpha)}
          {(p\wedge q_1)\wedge(v_{12}(p\wedge q_2)v_{21})\wedge\cdots\wedge
          (v_{1n}(p\wedge q_n)v_{n1})}, \tag13 $$
        so $h_1\le k$, $v_{1j}h_jv_{j1}\le k$ $(2\le j\le n)$.
      \enddemo

    We will compute $h_1\MvN h_1$ as a way of determining
      $q_1\MvN q_1=(q_1\NvN_2q_1\cup\{v_{1j}y_{j1}\mid2\le j\le n\})''$
      and hence $\MvN=(q_1\MvN q_1)\otimes M_n$.
    Let $z_j\in\{h_1,v_{1j}h_jv_{j1}\}''$ (the double commutant taken having
      $q_1$ for identity) be a partial isometry from $v_{1j}h_jv_{j1}$
      to $h_1$, $(2\le j\le n)$.
    Then $h_1\MvN h_1$ is generated by $h_1\NvN_2h_1$ together with
      $\{z_jv_{1j}y_{j1}\mid2\le j\le n\}$.

    \proclaim{Claim 3.2b} $z_jv_{1j}y_{j1}$ is a Haar unitary in $h_1\MvN h_1$
      for $2\le j\le n$.
    \endproclaim
    \Indeed
      Let $d=h_1-n(1-\alpha)q_1$, $c_j=v_{1j}\cw_jv_{j1}$ where
        $\cw_j=h_j-n(1-\alpha)q_j$ $(2\le j\le n)$.
      Then $z_j$ is the s.o.--limit of a bounded sequence in
        $\lspan\Lambda(\{d\},\{c_j\})$.
      Hence it suffices to show that a nontrivial alternating product of length
        $2m>0$ in $\Lambda(\{d\},\{c_j\})$ and $v_{1j}y_{j1}$ has trace zero.
      Noting that $\cw_j$ and $y_{j1}$ lie in $\NvN_1$, it suffices to show
        that an alternating product in $\{v_{1j},v_{j1}\}$ and $\NvN_1$
        that has $m$ more copies of $v_{1j}$ than of $v_{j1}$
        has trace zero.
      But we know that every element of $\NvN_1$ is the s.o.--limit of a
        bounded sequence in
        $\lspan\Lambda(\{a\},\{b_1,\ldots,b_n\})$.
      So it suffices to show that an alternating product in $\{v_{1j},v_{j1}\}$
        and $\Lambda(\{a\},\{b_1,\ldots,b_n\})$ that has $m$ more copies of
        $v_{1j}$ than of $v_{j1}$ has trace zero.
      But such a product can be reduced to give an
        alternating product in $\{a\}$ and
        $\{v_{1j},v_{j1},b_1,\ldots,b_n\}$, which must be nontrivial
        because there were more $v_{1j}$ than $v_{j1}$, so by freeness
        has trace zero.
      This proves the claim.
    \enddemo

    \proclaim{Claim 3.2c}
      $\bigl\{h_1\NvN_2h_1,(\{z_jv_{1j}y_{j1}\})_{2\le j\le n}\bigr\}$
        is $*$--free in $h_1\MvN h_1$.
    \endproclaim
    \Indeed
      Let $\NvN_{20}=(\{h_1\}\cup\{v_{1j}h_jv_{j1}\mid2\le j\le n\})''$
        (with identity element $q_1$).
      Let $Q=h_1\NvN_1h_1=L(\freeF_{n-1})$ and $R=h_1\NvN_{20}h_1$.
      Then $Q$ and $R$ together generate $h_1\NvN_2h_1$ (and by~1.2, $\{Q,R\}$
        is free in $h_1\NvN_2h_1$.
      We will show that $\{Q,R,(\{z_jv_{1j}y_{j1}\})_{2\le j\le n}\}$
        is free in $h_1\MvN h_1$.
      It suffices to show that a nontrivial traveling product in $\Ro$, $\Qo$,
        $(\{(z_jv_{1j}y_{j1})^m\mid m\in\Intg\backslash\{0\}\})_{2\le j\le n}$
        has trace zero.
      Regrouping gives a nontrivial alternating product, call it $\Pc$, in
        $$ \Omega_1=\Qo\cup\bigcup_{2\le j\le n}(Qy_{1j}\cup y_{j1}Q
          \cup y_{j1}\Qo y_{1j})\cup\bigcup\Sb2\le i,j\le n\\i\neq j\endSb
          y_{i1}Qy_{1j} $$
        and
        $$ \Omega_2=\Ro\cup\bigcup_{2\le j\le n}(v_{j1}z_j^*R\cup Rz_jv_{1j}
          \cup v_{j1}z_j^*\Ro z_jv_{1j})\cup\bigcup\Sb2\le i,j\le n\\
          i\neq j\endSb v_{j1}z_j^*R z_iv_{1i}, $$
        where an occurrence of $y_{1j}$ (be it from $Qy_{1j}$,
        $y_{j1}\Qo y_{1j}$ or $y_{i1}Qy_{1j}$) is always followed by an
        occurrence of $v_{j1}z_j^*$ (be it from $v_{j1}z_j^*R$,
        $v_{j1}z_j^*\Ro z_jv_{1j}$ or $v_{j1}z_j^*Rz_iv_{1i}$).
      We similarly have the list of conditions:
        $$ \matrix \format\c&\c&\c&\c&\c&\l\\
          \text{An occurrence of }&y_{1j}&\text{ is always }&\text{followed by}
            &\text{ an occurrence of }&v_{j1}z_j^*, \\
          \tq&y_{j1}&\tq&\tp&\tq&z_jv_{1j}, \\
          \tq&z_jv_{1j}&\tq&\tf&\tq&y_{j1}, \\
          \tq&v_{j1}z_j^*&\tq&\tp&\tq&y_{1j}.
        \endmatrix \tag14 $$
      Let $d,c_j,\cw_j$ be as in the proof of Claim~3.2b, and denote
        also $d=c_1$.
      Then $\{(c_j)_{1\le j\le n}\}$ is free in $h_1\MvN h_1$ (as we saw)
        and $\lspan\Lambda((\{c_j\})_{1\le j\le n})$ is a dense
        $*$--subalgebra of $\NvN_{20}$.
      For $1\le k,l\le n$, let $W_{kl}$ be the set of nontrivial traveling
        products in $(c_i)_{1\le i\le n}$ which begin with a letter other
        than $c_k$ and end with a letter other than $c_l$.
      If $x\in\Ro$, then x is the s.o.--limit of a bounded sequence in
        $h_1\lspan(W_{11})h_1$.
      And since $z_j\in\NvN_{20}$, an element $x\in Rz_j$ is the s.o.--limit
        of a bounded sequence in $h_1\lspan(\{1\}\cup W_{1j})v_{1j}h_jv_{j1}$,
        so every element of $Rz_jv_{1j}$ is the s.o.--limit of a bounded
        sequence in $h_1\lspan(\{1\}\cup W_{1j})v_{1j}h_j$.
      In a similar way we obtain approximating sequences for the other elements
        of $\Omega_2$, and hence deduce that it suffices to show that we
        get something of trace zero when we substitute in our alternating
        product $\Pc$
        $$ \matrix\format\c&\l&\l\\
          \text{for every occurrance of a letter from }&
          \Ro\text{ an arbitrary element of }&h_1W_{11}h_1, \\
          \tq&Rz_jv_{1j}\tqc&h_1(W_{1j}\cup\{1\})v_{1j}h_j, \\
          \tq&v_{j1}z_j^*R\tqc&h_jv_{j1}(W_{j1}\cup\{1\})h_1, \\
          \tq&v_{j1}z_j^*\Ro z_jv_{1j}\tqc&h_jv_{j1}(W_{jj})v_{1j}h_j \\
          \tq&v_{j1}z_j^*\Ro z_iv_{1i}\tqc
            &h_jv_{j1}(W_{ji}\cup\{1\})v_{1i}h_i\;(i\neq j). \\
        \endmatrix \tag15 $$
      When making the above substitutions, we do not need to worry about
        writing the $h_1$, $h_j$ or $h_i$ (from~(15)) because the neighboring
        element from $\Omega_1$ will absorb it, except possibly for an $h_1$
        at the very beginning and one at the very end.
      If an $h_1$ occurs at the beginning or the end, write
        $h_1=d+n(1-\alpha)q_1$ and distribute, giving a sum of up to four
        products.
      Now write each element of $W_{kl}$ ($1\le k,l\le n)$ as a word in
        $(\{c_i\})_{1\le i\le n}$, and for each $c_i$ with $i\ge2$ write
        $c_i=v_{1i}\cw_iv_{i1}$.
      But every element of $\Omega_1$ is the s.o.--limit of a bounded sequence
        in $\lspan(\Lambda(\{a\},\{b_1,\ldots,b_n\})
        \backslash\{1,b_1,\ldots,b_n\})$, as is $c_1$ and each
        $\cw_i$ $(2\le i\le n)$.
      So it suffices to show that replacing each letter from $\Omega_1$,
        each $c_1$ and each $\cw_i$ $(2\le i\le n)$ with an arbitrary
        element of
        $\Lambda(\{a\},\{b_1,\ldots,b_n\})\backslash\{1,b_1,\ldots,b_n\}$,
        we get something having trace zero.
      But because of the rules~(14) about neighbors, we see that by regrouping,
        we get a nontrivial alternating product in $\{a\}$ and
        $\{b_1,\ldots,b_n\}\cup\{v_{ij}\mid i\neq j\}$,
        which by freeness has trace zero.
      This proves Claim~3.2c.
    \enddemo

    From~(13) and Claims~3.2b and~c we can find $\MvN$ in the two cases
      to finish the proof of the proposition.
  \QED

  It may be of interest to note that the projection $p'$ in the statement
    of the proposition above is
    $$ p'=\sum_{1\le j\le n}(p\wedge q_j)\wedge
      \bigwedge\Sb 1\le i\le n\\i\neq j\endSb v_{ji}(p\wedge q_i)v_{ij}\le p.
$$

  Now we have all that we need to determine the free product of any given pair
    of finite dimensional algebras.
  It is most natural to phrase the results in terms of free dimension.

  \proclaim{Definition 3.3} \rm
    The {\it free dimension} of a finite dimensional algebra
      $$ A=\smd{M_{n_1}}{\alpha_1}\oplus\cdots\oplus\smd{M_{n_k}}{\alpha_k}
        \;\;(k\ge1) $$
      is equal to
      $$ \fdim(A)=\sum_{1\le i\le k}\alpha_i^2(1-n_i^{-2})
        +\sum\Sb1\le i,j\le k\\i\neq j\endSb\alpha_i\alpha_j. $$
    More generally, the {\it free dimension} of an algebra
      $$ B=\smd{L(\freeF_s)}{\alpha_0}\oplus\smd{M_{n_1}}{\alpha_1}
        \oplus\cdots\oplus\smd{M_{n_k}}{\alpha_k}
        \;\;(k\ge0,\,s\ge1) $$
      is equal to
      $$ \fdim(B)=\alpha_0^2s+\sum_{1\le i\le k}\alpha_i^2(1-n_i^{-2})
        +\sum\Sb0\le i,j\le k\\i\neq j\endSb\alpha_i\alpha_j. $$
  \endproclaim

  The arguments leading up tp Theorem~3.6 are once again essentially nothing
    more than induction resting on the shoulders of the results we have
    already proven.
  There are certainly many different ways to formally write down this
    procedure.
  The following sequence of lemmas outlines one such proof.
  Once again, we will omit the laborious and elementary algebraic
verifications.

  \proclaim{Lemma 3.4}
    Let
      $$ \Ac=\smd{L(\freeF_r)}{\alpha_0}\oplus\smd\Cpx{\alpha_1}
        \oplus\cdots\oplus\smdp\Cpx{\alpha_k}{p_k}
        \;\;(k\ge0,\,r\ge1,\,\alpha_0\ge0), $$
      where $0<\alpha_1\le\alpha_2\le\cdots\le\alpha_k$, and let $m\ge2$.
    Then
      $$ \Ac*M_m\cong\cases L(\freeF_s)&\text{if }\;\alpha_k\le1-m^{-2}\\
        L(\freeF_s)\oplus\smdp{M_m}{m^2\alpha_k-m^2+1}{p'}
        &\text{if }\;\alpha_k\ge1-m^{-2}, \endcases \tag16 $$
      where $s$ is such that the right hand side of~(16) has free dimension
      equal to $\fdim(\Ac)+1-m^{-2}$, and where $p'\le p_k$.
  \endproclaim
  \Pf
    Let $\MvN$ be the left hand side of~(16).
    The case $k=0$ follows from~3.1 and~2.4, as does the case $k\ge1$,
      $\alpha_k\le1-m^{-1}$.
    The case $m^{-2}\le\alpha_k$ follows from~1.2, 2.4 and~3.2 by writing
      $$ \matrix\format\c&\l&\c&\l\\
        \MvN&=(&\Bc&\oplus\smd\Cpx{\alpha_k})*M_m \\
        \cup \\
        \NvN&=(&\Cpx&\oplus\smd\Cpx{\alpha_k})*M_m
      \endmatrix $$
  \QED

  Similarly one can prove
  \proclaim{Proposition 3.5}
    Let
      $$ \Ac=\smd{L(\freeF_r)}{\alpha_0}\oplus\smdp{M_{n_1}}{\alpha_1}{p_1}
        \oplus\cdots\oplus\smdp{M_{n_k}}{\alpha_k}{p_k}
        \;\;(k\ge0,\,r\ge1,\,\alpha_0\ge0) $$
      and let $m\ge2$.
    If $n_j\ge2$ for all $1\le j\le k$, then
      $\Ac*M_m=L(\freeF(\fdim(\Ac)+1-m^{-2}))$.
    If some $n_j=1$, assume that $n_k=1$ and
      $\alpha_k=\max\{\alpha_j\mid1\le j\le k,\,n_j=1\}$.
    Then
      $$ \Ac*M_m=\cases L(\freeF(\fdim(\Ac)+1-m^{-2}))&\text{if }\;
        \alpha_k\le1-m^{-2} \\
        L(\freeF_s)\oplus\smdp{M_m}{m^2\alpha_k-m^2+1}{p'}
        &\text{if }\;\alpha_k\ge1-m^{-2}, \endcases \tag17 $$
      where $p'\le p_k$ and $s$ is such that the free dimension of the
      right hand side of~(17) is equal to $\fdim(\Ac)+1-m^{-2}$.
  \endproclaim

  \proclaim{Theorem 3.6}
    Let
      $$ \aligned
        &A=\smdp{M_{n_1}}{\alpha_1}{p_1}\oplus\cdots\oplus
          \smdp{M_{n_k}}{\alpha_k}{p_k} \\
        &B=\smdp{M_{m_1}}{\beta_1}{q_1}\oplus\cdots\oplus
          \smdp{M_{m_l}}{\beta_l}{q_l}
      \endaligned $$
      be finite dimensional algebras, each of dimension $\ge2$ and such that
      the sum of their dimensions is $\ge5$.
    Then
      $$ A*B=L(\freeF_s)\oplus\bigoplus\Sb1\le i\le k\\1\le j\le l\endSb
        \smdp M{\gamma_{ij}}{f_{ij}}_{N(i,j)}, \tag18 $$
      where $N(i,j)=\max(n_i,m_j)$,
      $\gamma_{ij}=N(i,j)^2\max(\dfrac{\alpha_i}{n_i^2}
      +\dfrac{\beta_j}{m_j^2}-1,0)$ and where $f_{ij}\le p_i\wedge q_j$.
    Note that $\gamma_{ij}>0$ implies either $n_i=1$ or $m_j=1$.
  \endproclaim
  \Pf
    By induction on $K=\card(\{n_i\ge2\}\cup\{m_j\ge2\})$.
    Let $\MvN$ be the left hand side of~(18).
    The case $K=0$ is just~2.3.
    For the inductive step, let $K\ge1$ and suppose $n_k\ge2$.
    Write $A=A_0\oplus\smd{M_{n_k}}{\alpha_k}$,
      $$ \matrix\format\c&\l&\c&\l\\
        \MvN&=(A_0\oplus&\smd{M_{n_k}}{\alpha_k}&)*B \\
        \cup \\
        \NvN&=(A_0\oplus&\smd\Cpx{\alpha_k}&)*B
      \endmatrix $$
      and use~1.2, the inductive hypothesis and~3.5.
  \QED

\noindent{\bf \S4. Hyperfinite algebras.}

  In this section we will extend the results of the last section
    to hyperfinite von Neumann algebras, {\it i.e\.} inductive limits
    of finite dimensional algebras (with trace--preserving morphisms).

  \proclaim{Definition 4.1} \rm
    Let $\Ac=L(\freeF_r)$, $\Bc=L(\freeF_{r'})$ with $r<r'$.
    Then $\psi:\Ac\rightarrow\Bc$ is a {\it standard embedding}
      if the following is satisfied:  for some W$^*$--noncommutative
      probability space $(\MvN,\phi)$ with $\phi$ a trace and with
      $\MvN$ containing a copy $R$ of the
      hyperfinite II$_1$--factor and a semicircular family
      $\omega=\{X^t\mid t\in T\}$ such that $\{R,\omega\}$ is free,
      there exist subsets $S\subset S'\subseteq T$, projections
      $p_s\in R$ $(s\in S')$ and isomorphisms
      $\alpha:\Ac\rightarrow(R\cup\{p_sX^sp_s\mid s\in S\})''$ and
      $\beta:\Bc\rightarrow(R\cup\{p_sX^sp_s\mid s\in S'\})''$ such that
      $\psi=\beta^{-1}\circ i\circ\alpha$, where $i$ is the inclusion.
  \endproclaim

  If $\psi:\Ac\rightarrow\Bc$ is an injective morphism of von Neumann
    algebras, $p\in\Ac$ a projection, then
    $\psi\restrict_{p\Ac p}:p\Ac p\rightarrow\psi(p)\Bc\psi(p)$ is also
    an injective morphism of von Neumann algebras.
  We have
  \proclaim{Proposition 4.2}
    Let $\Ac=L(\freeF_r)$, $\Bc=L(\freeF_{r'})$ with $r<r'$;
     suppose $\psi:\Ac\rightarrow\Bc$, $p\in\Ac$ a projection.
    Then $\psi$ is a standard embedding
      if and only if
      $\psi\restrict_{p\Ac p}$ is a standard embedding.
  \endproclaim
  \Pf
    Since tensoring with $M_n$
      preserves ``standardness,'' it is enough to prove one direction.
    Suppose $\psi$ is a standard embedding and let $\alpha$, $\beta$ be
      isomorphisms as in Definition~4.1.
    By the proof of~2.2 of~\cite{4}, we may assume that $p_s\le p$
      for all $s\in S'$.
    Then
      $$ \aligned
        \alpha\restrict_{p\Ac p}:&
          p\Ac p\isorarrow\bigl(pRp\cup\{p_sX^sp_s\mid s\in S\}\bigr)'' \\
        \beta\restrict_{\psi(p)\Bc\psi(p)}:&
          \psi(p)\Bc\psi(p)\isorarrow
          \bigl(pRp\cup\{p_sX^sp_s\mid s\in S'\}\bigr)''
      \endaligned \tag19 $$
      and $\psi\restrict_{p\Ac p}=(\beta\restrict_{\psi(p)\Bc\psi(p)})^{-1}
      i(\alpha\restrict_{p\Ac p})$.
    Theorem~1.3 of~\cite{4} tells us that the generators on the
      right hand side of~(19) are of the correct sort,
      so $\psi\restrict_{p\Ac p}$ is a standard embedding.
  \QED

  \proclaim{Proposition 4.3}
    \roster
    \item"(i)" The composition of standard embeddings is a standard embedding.
    \item"(ii)" Let $\Ac_n=L(\freeF_{r_n})$ for $n\ge1$ with $r_n<r_{n+1}$,
        $\psi_n:\Ac_n\rightarrow\Ac_{n+1}$  standard embeddings.
      Then $\dsize\varinjlim_n(\Ac_n,\psi_n)=L(\freeF_r)$
        where $r=\dsize\lim_{n\rightarrow\infty}r_n$.
    \endroster
  \endproclaim
  \Pf
    To prove~(i), let $\Ac=L(\freeF_r),$ $\Bc=L(\freeF_s)$ and
$\Cc=L(\freeF_t)$
      for $r<s<t$.
    Let $\psi:\Ac\rightarrow\Bc$ and $\psiw:\Bc\rightarrow\Cc$ be standard
      embeddings.
    Then taking $(\MvN,\phi)$, $\{X^t\mid t\in T\}$ as in Definition~4.1,
      there exist isomorphisms
      $$ \aligned
        \alpha:{}&\Ac\isorarrow(R\cup\{p_sX^sp_s\mid s\in S\})'' \\
        \beta:{}&\Bc\isorarrow(R\cup\{p_sX^sp_s\mid s\in S'\})'' \\
        \betaw:{}&\Bc\isorarrow(R\cup\{p_sX^sp_s\mid s\in \Sw\})'' \\
        \gammaw:{}&\Cc\isorarrow(R\cup\{p_sX^sp_s\mid s\in \Sw'\})'',
      \endaligned $$
      where $S\subset S'\subset T$, $\Sw\subset\Sw'\subset T$,
      $S'$ and $\Sw'$ are disjoint, $p_s$ are projections in $R$
      and $\psi=\beta^{-1}\circ i\circ\alpha$,
      $\psiw=\gammaw^{-1}\circ i\circ\betaw$.
    Using the isomorphism $\beta\circ\betaw^{-1}$ we get an isomorphism
      $$ \gammab:\Cc\isorarrow(R\cup\{p_sX^sp_s\mid s\in S'\}
        \cup\{\pw_sX^s\pw_s\mid s\in \SwpbSw\})'' $$
      where each $\pw_s\in(R\cup\{p_sX^sp_s\mid s\in S'\})''$
      and such that $\psiw\circ\psi=\gammab^{-1}\circ i\circ\alpha$.
    For $s\in\SwpbSw$ let $U_s\in(R\cup\{p_sX^sp_s\mid s\in S'\})''$
      be a unitary such that $U_s\pw_sU_s^*\in R$.
    We see that $\{R,(\{X^s\})_{s\in S'},(\{U_sX^sU_s^*\})_{s\in\SwpbSw}\}$
      is free and we may rewrite
      $$ \gammab:\Cc\isorarrow(R\cup\{p_sX^sp_s\mid s\in S'\}
        \cup\{(U_s\pw_sU_s^*)(U_sX^sU_s^*)(U_s\pw_sU_S^*)
        \mid s\in \SwpbSw\})'', $$
      thus proving~(i).

    In order to prove~(ii), we may apply the above argument recursively
      to obtain $S_1\subset S_2\subset S_3\subset\cdots$ and isomorphisms
      $$ \alpha_n:\Ac_n\isorarrow(R\cup\{p_sX^sp_s\mid s\in S_n\})'' $$
      such that $\psi_n=\alpha_{n+1}^{-1}\circ i\circ\alpha_n$.
    Then
      $$ \varinjlim_n(\Ac_n,\psi_n)
        \cong(R\cup\{p_sX^sp_s\mid s\in S_1\cup S_2\cup\cdots\})''
        =L(\freeF_r). $$
  \QED

  \proclaim{Proposition 4.4}
    The following are standard embeddings:
    \roster
    \item"(i)" the inclusion $i_1:L(\freeF_r)\hookrightarrow L(\freeF_r)
      *L(\freeF_{r'})$ for $1<r,r'\le\infty$;
    \item"(ii)" the inclusion $i_1:L(\freeF_r)\hookrightarrow L(\freeF_r)
      *B$ for $1<r\le\infty$, $B$ a finite dimensional algebra
      not equal to $\Cpx$;
    \endroster
  \endproclaim
  \Pf
    (i) In a noncommutative probability space $(\MvN,\phi)$ with $\phi$
      a trace, let $R$ and $\Rw$ be copies of the hyperfinite II$_1$--factor
      and $\omega=\{X^t\mid t\in T\}$ a semicircular family such that
      $\{R,\Rw,\omega\}$ is free.
    Let $S$, $\Sw$ be disjoint subsets of $T$ and $p_s\in R$ $(s\in S)$,
      $\pw_s\in\Rw$ $(s\in\Sw)$ be projections.
    Then it suffices to show that
      $$ \bigl(R\cup\{p_sX^sp_s\mid s\in S\}\bigr)''
        \overset i \to\hookrightarrow
        \bigl(R\cup\{p_sX^sp_s\mid s\in S\}\cup
        \Rw\cup\{\pw_sX^s\pw_s\mid s\in\Sw\}\bigr)'' \tag20 $$
      is a standard embedding.
    By~3.6 of~\cite{4}, there is a semicircular element $Y\in(R\cup\Rw)''$
      such that $\{R,\{Y\},\omega\}$ is free and $(R\cup\{Y\})''=(R\cup\Rw)''$.
    Then for each $s\in\Sw$ there is a unitary $U_s\in(R\cup\{Y\})''$
      such that $U_s\pw_sU_s^*=q_s\in R$.
    It is easily seen that
      $\{R,\{Y\},(\{X^s\})_{s\in S},(\{U_sX^sU_s^*\})_{s\in\Sw}\}$ is free.
    The right hand side of~(20) equals
      $(R\cup\{p_sX^sp_s\mid s\in S\}\cup\{Y\}\cup\{q_s(U_sX^sU_s^*)q_s
      \mid s\in\Sw\})''$, and the inclusion of~(20) is then clearly
      a standard embedding.

    (ii) Let $n\ge2$ be so large that $M_n*B$ is a factor and let
      $\{e_{ij}\mid1\le i,j\le n\}$ be a system of matrix units in
      $\Ac=L(\freeF_r)$.
    Then $\Ac=(e_{11}\Ac e_{11})\otimes M_n$.
    So
      $$ i_1\restrict_{e_{11}\Ac e_{11}}: e_{11}\Ac e_{11}\hookrightarrow
        e_{11}(\Ac*B)e_{11} $$
      is by~1.2 conjugate to the inclusion
      $$ e_{11}\Ac e_{11}\hookrightarrow(e_{11}\Ac e_{11})
        *(e_{11}(M_n*B)e_{11}). \tag21 $$
    But by the compression formula~(2) we have
      $e_{11}\Ac e_{11}=L(\freeF(1+n^2(r-1)))$,
      and by~3.6 together with the compression formula we have
      $e_{11}(M_n*B)e_{11}=L(\freeF(n^2\fdim(B)))$.
    Then by~(i),~(21) is a standard embedding, so by~4.2,
      $\Ac\hookrightarrow\Ac*B$ is a standard embedding.
  \QED

  Let $\Ac$ be a hyperfinite von Neumann algebra (with implicitly specified
    trace).
  Then decomposing $\Ac$ over its center, we may write
    $$ \Ac=\smd\Acw\alpha\oplus
     \bigoplus_{i\in I}\smdp{M_{n_i}}{\alpha_i}{p_i}, \tag22 $$
    where $\alpha\ge0$, $I$ is a finite or countably infinite (or empty)
    index set, $n_i\in\Nats\backslash\{0\}$ and $\alpha_i>0$,
    and where $\Acw$ is a diffuse von Neumann algebra, which can thus be
    further written as
    $$ \Acw=\smd{L(\Intg)}{\gamma_1}\oplus\smd{(L(\Intg)\otimes M_2)}{\gamma_2}
      \oplus\smd{(L(\Intg)\otimes M_3)}{\gamma_3}\oplus\cdots\oplus
      \smd{(L(\Intg)\otimes R)}{\gamma_\infty}\oplus
      \bigoplus_{j\in J}\smd R{\delta_j}, $$
    where $R$ is the hyperfinite II$_1$--factor,
    $\gamma_k\ge0$, $J$ is a finite
    or countably infinite (or empty) index set and $\delta_j>0$.

  \proclaim{Definition 4.5}\rm
    The {\it free dimension} of a hyperfinite von Neumann algebra $\Ac$ written
      as in~(22) is equal to
      $$ \fdim(\Ac)=\alpha^2+\sum_{i\in I}\alpha_i^2(1-n_i^{-2})
        +2\alpha(1-\alpha)+\sum\Sb i,j\in I\\i\neq j\endSb\alpha_i\alpha_j. $$
  \endproclaim

  \proclaim{Theorem 4.6}
    Let $\Ac$ and $\Bc$ be hyperfinite (or finite dimensional)
      von Neumann algebras such that $\dim(\Ac),\dim(\Bc)\ge2$ and
      $\dim(\Ac)+\dim(\Bc)\ge 5$.
    Let $\Ac$ be written as in~(22) and write
      $$ \Bc=\smd\Bcw\beta\oplus
        \bigoplus_{j\in J}\smdp{M_{m_j}}{\beta_j}{q_j}, $$
      where $\beta\ge0$, $J$ is a finite or countably infinite (or empty)
      index set, $m_j\in\Nats\backslash\{0\}$ and $\beta_j>0$,
      and where $\Bcw$ is a diffuse von Neumann algebra.
    Then
      $$ \Ac*\Bc\cong \smd{L(\freeF_s)}\gamma\oplus
        \bigoplus\Sb i\in I\\j\in
J\endSb\smdp{M_{N(i,j)}}{\gamma_{ij}}{f_{ij}},
        \tag23 $$
      where $N(i,j)=\max(n_i,m_j)$,
      $\gamma_{ij}=N(i,j)^2\max(\dfrac{\alpha_i}{n_i^2}
      +\dfrac{\beta_j}{m_j^2}-1,0)$
      and $f_{ij}\le p_i\wedge q_j$, and where $s$ is such that the free
      dimension of the right hand side of~(23) equals $\fdim(\Ac)+\fdim(\Bc)$.
    Note that $\gamma_{ij}>0$ implies either $n_i=1$ or $m_i=1$, and that
      there are only finitely many pairs $(i,j)\in I\times J$ for which
      $\gamma_{ij}>0$.
  \endproclaim
  \Pf
    Let $(A_k)_{k\ge1}$ be an increasing sequence of finite dimensional
      subalgebras of $\Ac$ such that $\cup A_k$ is dense in $\Ac$ and
      such that each inclusion $A_k\hookrightarrow A_{k+1}$ (which we
      denote $\phi_k$) is of the form
      $A_k=M_n\oplus A\hookrightarrow(M_n\otimes C)\oplus A=A_{k+1}$, where
      $A$ and $C$ are finite dimensional algebras (depending on $k$).
    Then $\fdim(A_k)$ is a non--decreasing sequence whose limit is
$\fdim(\Ac)$.
    (If $\Ac$ is finite dimensional then eventually the $A_k$ are all the
same.)
    Let $(B_k)_{k\ge1}$ be a similar sequence for $\Bc$, and denote the
      inclusion $B_k\hookrightarrow B_{k+1}$ by $\psi_k$.
    Theorem~3.6 allows us to find each $A_k*B_k$, and for large enough $k$, we
      have that
      $$ A_k*B_k=\smd{L(\freeF_{s_k})}\gamma\oplus
        \bigoplus\Sb i\in I\\j\in J\endSb
        \smdp{M_{N(ij)}}{\gamma_{ij}}{f_{ij}}, $$
      with $N(i,j)$ and $\gamma_{ij}$ as in~(23).
    Moreover, letting $f=1-\sum_{i\in I,j\in J}f_{ij}$, we see that
      $\phi_k*\psi_k\restrict_{(1-f)(A_k*B_k)(1-f)}$ is the identity
      map on the finite dimensional part $\bigoplus M_{N(ij)}$.
    Using~1.4, 4.2, 4.3(i) and~4.4(ii), we can show that
      $\phi_k*\psi_k\restrict_{f(A_k*B_k)f}:L(\freeF_{s_k})
      \hookrightarrow L(\freeF_{s_{k+1}})$ is a standard embedding.
    Taking inductive limits and applying~4.3(ii) proves the theorem.
  \QED

  \proclaim{Remark 4.7}\rm
    The precise condition for $\Ac*\Bc$ to be a factor may be described
      as follows.
    Let $\Ac$ be as written in~(22) and define the {\it lumpiness} of $\Ac$
      to be $l_\Ac=\max\{\tfrac{\alpha_i}{n_i^2}\mid i\in I\}$.
    Let $l_\Bc$ be the lumpiness of $\Bc$.
    Then $\Ac*\Bc$ is a factor if and only if $l_\Ac+l_\Bc\le1$,
      (provided that neither $\Ac$ nor $\Bc$ is one--dimensional and one
      of them is $\ge3$--dimensional).
  \endproclaim

\noindent{\bf \S5. Group algebras.}

  In the fundamental paper~\cite{2},
    A\. Connes showed that $L(G)$ is hyperfinite for $G$
    a discrete amenable group.
  It is now easy to find the free dimension of  such $L(G)$, and
    to that end we want to have the following proposition which has been
    known for some time.
    \footnote{We would like to thank Sorin Popa:  when asked if this result
    was true, he communicated to us that he had proved it in 1981, and he
    made helpful comments which enabled us to find this proof.}

  \proclaim{Proposition 5.1}
    Let $G$ be an infinite discrete group.
    Then $L(G)$ is diffuse.
  \endproclaim
  \Pf
    Suppose for contradiction that
      $p\in L(G)$ is a central projection having trace $\alpha>0$
      and that $pL(G)\cong M_n$.
    Then $p(l^2(G))\cong L^2(pL(G),\alpha^{-1}\tau)\cong L^2(M_n,\tau_2)$,
      (where $\tau$ is the canonical trace on $L(G)$, $\tau_2$ the
      normalized trace on $M_n$),
      so $p(l^2(G))$ has dimension $n^2$.
    Let $U_g$ for $g\in G$ denote the left translation operator on $l^2(G)$
      and let $v_1,\ldots,v_{n^2}$ be an orthonormal basis for $p(l^2(G))$.
    Then since $v_i,v_j\in l^2(G)$ we have
      $\langle U_g v_i,v_j\rangle\rightarrow0$ as $g\rightarrow\infty$.
    But this contradicts the fact that $U_g\restrict_{p(l^2(G))}$
      is unitary for all $g\in G$.
  \QED

  \proclaim{Proposition 5.2}
    Let $G$ be an amenable discrete group.
    Then $L(G)$ has free dimension equal to $1-|G|^{-1}$
      (where $\infty^{-1}=0$).
  \endproclaim
  \Pf
    For $|G|=\infty$, we have by~5.1 that $L(G)$ is a diffuse hyperfinite
      von Neumann algebra, so has free dimension equal to $1$.
    For $|G|<\infty$, it is well-known that
      $$ L(G)=\bigoplus_{\alpha\in\Ghat}\smd{M_{n_\alpha}}{\gamma_\alpha}, $$
      where $\Ghat$ is the collection of equivalence classes of irreducible
      representations of $G$, $n_\alpha$ is the dimension of the representation
      and the tracial weight $\gamma_\alpha$ equals $n_\alpha^2/|G|$.
    One then uses~3.3.
  \QED

  \proclaim{Corollary 5.3}
    Let $G$ and $H$ be discrete amenable groups.
    Then $L(G*H)$ depends only on the orders of $G$ and $H$.
    More precisely,
      $$ L(G*H)=L(\freeF(2-|G|^{-1}-|H|^{-1}))
        \;\;(\text{where }\infty^{-1}=0) $$
      provided $|G|,|H|\ge2$ and $|G|+|H|\ge5$.
  \endproclaim

  \proclaim{Corollary 5.4}
     \footnote{Thanks to Florin Boca for pointing this out after reading
        a preliminary version of this paper.}
    Let $G_n$ ($n=1,2,\ldots$) be nontrivial discrete amenable (or finite)
    groups.
    Then
      $$ L(\operatornamewithlimits{*}_{n=1}^\infty G_n)=L(\freeF_\infty). $$
  \endproclaim
  \Pf
    $\Ac=L({\dsize\operatornamewithlimits * _{n=1}^\infty}G_n)=
      {\dsize\operatornamewithlimits * _{k=1}^\infty}L(G_{2k-1}*G_{2k})$.
    Since $\Intg_2*\Intg_2$ is amenable, we may assume without loss
      of generality that for no $k$ are both
      $G_{2k-1}$ and $G_k$ of order two.
    Then $G_{2k-1}*G_{2k}=
      L(\freeF(2-|G_{2k-1}|^{-1}-|G_{2k}|^{-1})$ and
      $L(\freeF(2-|G_{2k-1}|^{-1}-|G_{2k}|^{-1})=(R_k\cup\{p_kX_kp_x\})''$
      where $R_k$ is a copy of the hyperfinite II$_1$--factor, $p_k\in R_k$ is
      a projection, $X_k$ is a semicircular element
      and $\{R_k,\{X_k\}\}$ is free.
    Thus $\Ac\cong(\bigcup_{k=1}^\infty R_k\cup\{p_kX_kp_k\mid k\ge1\})''$
      where $R_k$, $p_k$ and $X_k$ are as above and
      $\{(R_k)_{k\ge1},(\{X_k\})_{k\ge1}\}$ is free.
    Now $\Ac_0=(\bigcup_{k\ge1}R_k)''\cong
      {\dsize\operatornamewithlimits * _{k=1}^\infty}R_k=
      {\dsize\operatornamewithlimits * _{j=1}^\infty}(R_{2j-1}*R_{2j})
      \cong{\dsize\operatornamewithlimits * _{j=1}^\infty}L(\freeF_2)=
      L(\freeF_\infty)=(R\cup\{Y_l\mid l\ge1\})''$
      where $R$ is a copy of the hyperfinite II$_1$--factor,
      $Y_l$ are semicircular elements and $\{R,(\{Y_l\})_{l\ge1}\}$ is free.
    For each $k$ there is a unitary $U_k\in\Ac_0$ such that $U_kp_ku_k^*\in R$.
    Then $\Ac\cong(R\cup\{Y_l\mid l\ge1\}\cup
      \{(U_kp_kU_k^*)(U_kX_kU_k^*)(U_kp_kU_k^*)\mid k\ge1\})''$
      and it is easily seen that
      $\{R,(\{Y_l\})_{l\ge1},(\{U_kX_kU_k^*\})_{k\ge1}\}$ is free, so
      (by Defenition~2.1 of~\cite{4}) $\Ac=L(\freeF_\infty)$.
  \QED

\smallpagebreak
\noindent{\bf Acknowledgements.}

  I would like to thank Dan Voiculescu, my advisor, for helpful comments,
    also Dietmar Bisch for helpful conversations, Florin Boca for
    conversations about free products of finite dimensional abelian algebras
    and his (correct) conjecture about when they are factors
    and Sorin Popa for helpful comments about group algebras.

\end{document}